\documentclass[12pt]{article}
\usepackage{a4wide,amsmath,amssymb,hyperref}
\newcommand{\ra}{R_{\alpha}}
\newcommand{\La}{\Lambda}
\newcommand{\De}{\Delta}

\newcommand{\p}{\psi}

\newtheorem{theorem}{Theorem}
\newtheorem{conjecture}{Conjecture}

\newtheorem{corollary}[theorem]{Corollary}

\newtheorem{duffinschaeffer}{The Duffin-Schaeffer conjecture: }

\newtheorem{catlin}{Catlin's conjecture: }

\newtheorem{catlinsim}{The simultaneous Catlin conjecture: }

\newtheorem{catlinlin}{The linear forms Catlin conjecture: }

\newtheorem{wald}{Waldschmidt's conjecture: }

\newtheorem{DSdual}{The dual Duffin-Schaeffer conjecture: }

\newtheorem{DSsystems}{The linear forms Duffin-Schaeffer conjecture: }

\renewcommand{\Bbb}[1]{\mathbb{#1}}
\newcommand{\A}{{\Bbb A}}         
\newcommand{\I}{{\Bbb I}}         
\newcommand{\N}{{\Bbb N}}         
\newcommand{\Q}{{\Bbb Q}}         
\newcommand{\R}{{\Bbb R}}         
\newcommand{\Rp}{\R^+}    
\newcommand{\Z}{{\Bbb Z}}         

\newcommand{\ds}{\displaystyle}

\def\SBA{\cB}
\def\SWA{\cW}
\def\DS{Duffin-Schaeffer}

\def\Ja{Jarn{\'\i}k}
\def\JB{Jarn{\'\i}k-Besicovitch}
\def\K{Khintchine}

\def\Jarnik{Jarn\'\i k}
\def\MTP{Mass Transference Principle}
\newcommand{\cA}{{\cal A}}
\newcommand{\cB}{{\cal B}}

\newcommand{\cD}{{\cal D}}

\newcommand{\cF}{{\cal F}}

\newcommand{\cH}{{\cal H}}

\newcommand{\cN}{{\cal N}}

\newcommand{\cR}{{\cal R}}

\newcommand{\cV}{{\cal V}}
\newcommand{\cW}{{\cal W}}

\newcommand{\cZ}{{\cal Z}}

\newcommand{\ve}{\varepsilon}

\newcommand{\diam}{\operatorname{diam}}
\newcommand{\dist}{\operatorname{dist}}

\newcommand{\vv}[1]{{\mathbf{#1}}}

\renewcommand{\le}{\leqslant}
\renewcommand{\ge}{\geqslant}
\newcommand{\nz}{\smallsetminus\{0\}}

\def\DA{Diophantine approximation}
\def\ie{\emph{i.e.}}

\def\xb{\mathbf{x}}

\def\pb{\mathbf{p}}
\def\q{\mathbf{q}}

\begin{document}

\title{Classical metric Diophantine approximation revisited}

\author{Victor Beresnevich\footnote{EPSRC Advanced Research Fellow, EP/C54076X/1}
\\ {\small\sc York} \and
Vasily Bernik
\\ {\small\sc Minsk} \and
Maurice Dodson \\ {\small\sc York} \and Sanju Velani\\
{\small\sc York}}

\date{{\small\it Dedicated to Klaus Roth on the occasion of his 80th birthday}}

\maketitle

\begin{abstract} The idea of using measure theoretic concepts
  to investigate the size of number theoretic sets, originating with
  E. Borel, has been used for nearly a century.  It has led to the
  development of the theory of metrical Diophantine approximation, a
  branch of Number Theory which draws on a rich and broad variety of
  mathematics.
  We discuss some recent progress and open problems concerning this classical
  theory. In particular, generalisations of the Duffin-Schaeffer and Catlin conjectures
  are formulated and explored.
\end{abstract}

\section{Dirichlet, Roth and the metrical
theory}\label{sec:Dirichl-Roth-metrica-theory}

Diophantine approximation is based on a quantitative analysis of the
property that the rational numbers are dense in the real line.
Dirichlet's theorem, a fundamental result of this theory, says that
given any real number $x$ and any natural number $N$, there are
integers $p$ and $q$ such that
$$
\left|x-\frac pq\right|\le\frac 1{qN}, \qquad 0<q< N\,.
$$
There is an extraordinarily rich variety of analogues and
generalisations of this fact. Although simple, Dirichlet's result is
best possible for all real numbers. Also it implies that for any
irrational $x$ there are infinitely many rational numbers $p/q$
satisfying $ \left|x-p/q\right|<1/{q^2}\,. $ The latter inequality can
be sharpened a little by the multiplicative factor $1/\sqrt 5$ but no
further improvement is possible for the golden ratio $(\sqrt5-1)/2$
and its equivalents \cite{Casselshort}. This leads to the notion of badly
approximable numbers: the irrational number $x$ is \emph{badly
  approximable}\/ if there is a constant $c>0$ such that
\begin{equation}\label{e:001}
    \left| x -\frac pq\right|<\frac c{q^2}
\end{equation}
for no rationals $p/q$. Dirichlet's theorem can be improved for
and only for badly approximable numbers~\cite{DavSchmidt70a}.  In
the
  other direction, if for any constant $c>0$ there is a rational $p/q$
  satisfying (\ref{e:001}) then $\beta$ is called \emph{well
    approximable}. The sets of badly and well approximable numbers
  will be denoted by $\SBA$ and $\SWA$ respectively.  Further,
  \emph{very well approximable} numbers $x$ satisfy the stronger
  condition that there exists an $\ve=\ve(x) > 0 $ such that the
  inequality
  $$
    \left|\alpha-\frac pq\right|<\frac 1{q^{2+\ve}}
  $$
holds for infinitely many rationals $p/q$. The set of such numbers
will be denoted by $\cV$. The numbers which are not very well
approximable will be called \emph{relatively badly approximable}\/
and will be denoted by $\mathcal{R}$, thus $\cR:=\R\setminus\cV$.
One readily sees that $\SBA\subset\mathcal{R}$.  In 1912
Borel~\cite{Borel12} proved that $\cV$ is of Lebesgue measure zero
(or null), so that its complement $\cR$ is of full measure.  The
Lebesgue measure of set $A$ will be written as $|A|$; there should
be no confusion with the symbol for modulus. As usually, we say
that `almost no' number lies in $\cV$ and `almost all' numbers lie
in its complement $\cR$ or that numbers in $\cR$ are `typical'.
If numbers in a unit interval are considered, there is a natural
interpretation in terms of probability: a number lies in $\cV$
with probability zero and in $\cR$ with probability one.

There are uncountably many badly approximable numbers, as they are
characterized by having bounded partial quotients in their
continued fraction expansion. This characterisation also implies
that $\SBA$ has Lebesgue measure zero (this also follows from
Khintchine's theorem discussed below) and full Hausdorff
dimension: $\dim\cB=1$. The quadratic irrationals are known to be
badly approximable and a natural but as yet unanswered question is
whether other irrational algebraic numbers \emph{are}\/ or
\emph{are not}\/ badly approximable. $\mathcal{R}\!$\textrm{oth}
proved the following remarkable result.
\begin{theorem}[Roth]
All real irrational algebraic numbers are $\mathcal{R}$-numbers.
\end{theorem}
From the metrical point of view, Roth's theorem ~\cite{Roth} shows
that every real algebraic irrational number behaves typically.

\section{Khintchine's theorem}\label{K}

Loosely speaking, in the above section we have been dealing with
variations of Dirichlet's theorem in which the right hand side or error of
approximation is either of the form $cq^{-2}$ or
$q^{-2-\ve}$. It is natural to broaden the discussion to general
error functions. More precisely,
given a function  $\psi:\N\to\Rp$, a real number $x$ is said
to be \emph{$\psi$--approximable} if there are infinitely many
$q\in\N$ such that
\begin{equation}\label{e:002}
    \|qx\|<\psi(q)  \ .
\end{equation}
The function $\psi$ governs the `rate' at which the rationals
approximate the reals and will be referred to as an
\emph{approximating function}. Here and throughout, $\|x\|$ denotes the
distance of $x$ from the nearest integer and $\Rp=[0,\infty)$.
One can readily verify
that the set of $\p$-approximable numbers is invariant under
translations by integer vectors. Therefore without any loss of
generality and to ease the `metrical' discussion which follows, we
shall restrict our attention to  $\psi$--approximable numbers in the
unit interval $\I:=[0,1)$.  The set of such numbers is clearly a
subset of $\I$ and  will be denoted by $\cA(\psi)$; i.e.
  $$
  \cA(\psi): =
\{x\in \I \colon \|qx\|<\psi(q) \text{ for infinitely many }
q\in\N \}  \ .
$$
In 1924, \K~\cite{Kh24} established a beautiful and strikingly
simple criterion for the `size' of the set $\cA(\psi)$ expressed
in terms of Lebesgue measure.  We will write the $n$-dimensional
Lebesgue measure of a set $X$ in $\R^n$ by $|X|_n$; when there is
no risk of confusion the suffix will be omitted.  There should be
no confusion with the notation for a norm or modulus of a number
or a vector. We give an improved modern version of this
fundamental result -- see
\cite{Casselshort,HarmanMNT,Schmidtbook,Sprindzuk}.

\begin{theorem}[Khintchine]
$$
    |\cA(\psi)| \ = \ \begin{cases} 0
      &\ds\text{if } \quad \sum_{r=1}^\infty \, \psi(r)<\infty \ , \\[3ex]
      1
      &\ds\text{if } \quad  \sum_{r=1}^\infty \, \psi(r)=\infty
      \  \text{ and $\psi$ is
                    monotonic}.
                  \end{cases}
$$
\end{theorem}

\medskip

\textit{Remark.\ } Regarding the above theorem and indeed the theorems and conjectures below,
it is straightforward to establish the complementary convergent statements;
i.e. if the sum in question converges then the set in question is of zero measure.

\medskip

In Khintchine's theorem, the divergence case constitutes the main
substance and involves the extra monotonicity condition. This
condition cannot in general be relaxed, as was shown by Duffin and
Schaeffer \cite{DuffinSchaeffer} in 1941. They constructed a
non-monotonic approximating function $\vartheta$ for which the sum
$\sum_q \vartheta(q)$ diverges but $|\cA(\vartheta)|_1=0$
(see~\cite{HarmanMNT,Sprindzuk} for details). Nevertheless, in the
case of arbitrary $\psi$, Duffin and Schaeffer produced a
conjecture that we now discuss.

The integer $p$ implicit in the inequality (\ref{e:002}) satisfies
\begin{equation}\label{e:003}
\left|x-\frac pq\right|<\frac{\psi(q)}q\,.
\end{equation}
To relate the rational $p/q$ with the error of approximation
$\psi(q)/q$ uniquely, we impose the coprimeness condition
$(p,q)=1$. In this case, let  $\cA'(\psi)$ denote the set of $x$
in $\I$ for which the inequality (\ref{e:003}) holds for
infinitely many $(p,q)\in\Z\times\N$ with $(p,q)=1$.  Clearly,
$$\cA'(\psi)\subset\cA(\psi)$$ so that the convergence part of
Khintchine's theorem remains valid if $\cA(\psi)$ is replaced by
$\cA'(\psi)$. In fact, for any approximating function $\psi\colon
\N\to \Rp$ one easily deduces  that
$$
|\cA'(\psi)|=   0  \quad {\rm if}  \quad
\displaystyle\sum_{r=1}^\infty \, \varphi(r)  \ \dfrac{\psi(r)}r \
< \ \infty \ .
$$
Here, and throughout, $\varphi$ is the Euler function.  A less
obvious fact is that the divergence part of Khintchine's theorem
(when $\psi$ is required to be monotonic) holds for $\cA'(\psi)$,
\ie, for $\psi$ monotonic, the coprimeness condition $(p,q)=1$ is
irrelevant. As already mentioned above, this is not the case if we
remove the monotonicity condition and the appropriate statement is
given by a famous conjecture.

\begin{duffinschaeffer}
For any approximating function $\psi\colon \N\to \Rp$
$$
|\cA'(\psi)|=   1  \quad if \quad
\displaystyle\sum_{r=1}^\infty\varphi(r) \dfrac{\psi(r)}r=\infty \
.
$$
\end{duffinschaeffer}

\noindent Although various partial results have been established
-- see \cite{HarmanMNT} for  details and references, the full
conjecture represents one of the most difficult and profound
unsolved problems in metric number theory.

We now turn our attention to the `raw' set $\cA(\psi)$ on which no
monotonicity or coprimeness conditions are imposed. It is known (see
\cite{HarmanMNT,Sprindzuk}) that the one-dimensional Lebesgue
measure of $\cA(\psi)$ is either zero or one  but this is far from
providing a criterion for $|\cA(\psi)|$ akin to the Duffin-Schaeffer
conjecture for $\cA'(\psi)$ or Khintchine's theorem for $\cA(\psi)$
with $\psi$ monotonic.
The following conjecture \cite{catlin76} provides such a criterion.

\begin{catlin}
For any approximating function $\psi\colon \N\to \Rp$
$$
|\cA(\psi)|=   1  \quad if \quad
\displaystyle\sum_{r=1}^\infty  \varphi(r)  \;
\max_{t\ge1} \dfrac{\psi(rt)}{rt} =\infty \ .
$$
\end{catlin}

\noindent Thinking geometrically, given a rational point $s/r$ with
$(s,r) = 1$ consider all its representations $p/q$ with $p=ts$ and
$q=tr$ for some $t\in\N$. The length of the largest interval given
by (\ref{e:003}) is precisely the maximum term appearing in the
above conjecture.

To the best of our knowledge, it is not known whether the
Duffin-Schaeffer conjecture is equivalent to  Catlin's conjecture.
For the current situation regarding the conjectures of
Duffin $\&$ Schaeffer and Catlin see~\cite[pp.~27--29]{HarmanMNT}.
It is remarkable that in the simultaneous setup
the analogues of the these conjectures have been
completely settled -- see \S\ref{sim}.

\subsection{Back to Roth}

We end this section by returning to Roth. In connection with Roth's
theorem, Waldschmidt \cite[pp.~260]{Waldschmidt} has made the
following conjecture which currently seems well beyond reach.

\begin{wald}
Let $\psi$ be a monotonic approximating  function  such that
$\sum_{r=1}^{\infty} \psi(r)<\infty $ and let $\alpha \in \I$ be a
real algebraic irrational number.  Then, $\alpha\not\in\cA(\psi)$.
\end{wald}

\noindent The conjecture is a general version of  Lang's
conjecture \cite{Lang}. The latter corresponds to the case that
$\psi : r \to r^{-1} (\log q)^{-1-\ve}$ with $\ve>0 \, $
arbitrary; i.e. the inequality
$$
\|q\alpha\| \, <  \, q^{-1}  ( \log q)^{-1-\ve}
$$
has only finitely many solution for every positive $\ve$. Note
that in view of the imposed convergent sum condition, we have that
$|\cA(\psi)|=0$. Thus, from the metrical point of view,
Waldschmidt's conjecture simply states that $\alpha$ behaves
typically in the sense that $\alpha$ belongs to the set
$\I\setminus\cA(\psi)$ of full measure. This clearly strengthens
the notion of typical as implied by Roth's theorem.

Within the statement of Waldschmidt's conjecture, it is natural to question the relevance
of the monotonicity assumption. In other words, does it make sense to consider the
following stronger form of the conjecture?  {\em  Let $\psi$ be an  approximating function such that
$\sum_{r=1}^{\infty} \varphi(r)  \,  \psi(r) / r <\infty $ and let
$\alpha \in \I$ be a real algebraic irrational number. Then,
$\alpha\not\in\cA(\psi)$.} This statement is easily seen to be false.
 For $q \in \N$, let
$\psi^*(q):=\|q\alpha\|$. Obviously,
$\liminf_{q\to\infty}\psi^*(q)=0$. Therefore, there exists a
sequence $\{q_n\}_{n\in\N}$ such that
$\sum_{n=1}^\infty\psi^*(q_n)<\infty$. Now, set $\psi(q)$ to be
$2\psi^*(q)$ if $q=q_n$ for some $n$ and $0$ otherwise. Thus for the non-monotonic approximating function $\psi$, we have that
$$\textstyle{\sum_{r=1}^{\infty} \varphi(r)  \,  \psi(r) / r  \ < \  \sum_{r=1}^\infty\psi(r)<\infty  \quad {\rm  but } \quad   \alpha \in\cA(\psi)  \ . } $$

\noindent The upshot of this is that the monotonicity assumption in Waldschmidt's conjecture
can not be removed. However, the example given above is somewhat `artificial' and it makes perfect sense
to study `density' questions of the following type: {\it for a given
approximating function $\psi$ how large is the set of
$\psi$-approximable algebraic numbers of degree $n$ compared to the
set of all algebraic numbers of degree $n$?} Even for
monotonic approximating functions, considering this question would
give a partial `metrical' or `density' answer to Waldschmidt's
conjecture.
%
%
%
%
%

%


\section{Simultaneous approximation by rationals}\label{sim}

In simultaneous \DA, one considers the set $\cA_m(\psi)$ of points
$\xb =(x_1,\dots,x_m)\in \I^m := [0,1)^m$ for which the inequality
 $$
 \|q\xb\|:=\max\left\{||qx_1||,\dots,||qx_m|| \right\} < \psi(q)
 $$
 holds for infinitely many positive integers $q$.   \K{} extended his
one-dimensional result discussed in \S\ref{K} to simultaneous
approximation~\cite{Kh25} (see
also~\cite[Chapter~VII]{Casselshort}).

\begin{theorem}[Khintchine]
$$
    |\cA_{m}(\psi)|_m = \begin{cases} 0
      &\ds\text{if } \ \sum_{r=1}^\infty \psi(r)^m<\infty \ , \\[3ex]
      1
      &\ds\text{if } \  \sum_{r=1}^\infty \psi(r)^m=\infty
      \   \text{ and $\psi$ is
                    monotonic}.
                  \end{cases}
$$
\end{theorem}

\noindent As with the one-dimensional statement,  Khintchine
originally had a stronger monotonicity condition. It turns out that
the monotonicity condition can be safely removed from Khintchine's
theorem for $m\ge2$. In fact, that this is the case, is a simple
consequence of the next result that deals with the setup in which a
natural coprimeness condition on the rational approximates is
imposed.

Let $\cA_m'(\psi)$ denote the set of points $\xb := (x_1,\dots,x_m) \in \I^m$ for
which the inequality
\begin{equation}\label{e:004}
 \left|\xb-\frac{\pb}{q}\right| < \frac{\psi(q)}{q} \quad {\rm with } \quad  (p_1,\dots,p_m,q)=1
\end{equation}
 is satisfied  for infinitely many $(\vv
p,q)\in\Z^m\times\N$.
The coprimeness condition imposed in the definition of
$\cA_m'(\psi)$ ensures that the points in $\I^m$ are approximated by
distinct rationals; i.e. the points $ \pb/q   :=(p_1/q,
\ldots,p_m/q)$ are distinct. The following metric result concerning
the set $\cA_m'(\psi)$ is due to Gallagher \cite{Gallagher65} and is
free from any monotonicity condition.

\begin{theorem}[Gallagher]
Let $m\ge2$. Then
$$
    |\cA'_{m}(\psi)|_m = \begin{cases} 0
      &\ds\text{if } \ \sum_{r=1}^\infty \psi(r)^m<\infty  , \\[3ex]
      1
      &\ds\text{if }  \ \sum_{r=1}^\infty \psi(r)^m=\infty .
                  \end{cases}
$$
\end{theorem}

To see that Gallagher's theorem removes the monotonicity
requirement from Khintchine's theorem for $m\ge2$, simply note
that divergent/convergent sum condition is the same in both
statements and that $ \cA'_{m}(\psi) \subset \cA_{m}(\psi) $.  The
latter implies that  $ |\cA_{m}(\psi)|_m = 1$ whenever $
|\cA'_{m}(\psi)|_m = 1$. Thus, for $m \ge 2$ we are able to
establish a criterion for the size of the `raw' set
$\cA_{m}(\psi)$ on which no monotonicity or coprimeness conditions
are imposed. In particular, we have the following divergent
statement as a corollary to Gallagher's theorem.

\begin{corollary}
\label{corgal}
Let $m\ge2$.  For any approximating function $\psi\colon \N\to \Rp$
$$
|\cA_m(\psi)|_m=   1  \quad if \quad
\displaystyle \sum_{r=1}^\infty \psi(r)^m=\infty \ .
$$
\end{corollary}

\noindent This  naturally  settles the
(simultaneous) higher dimensional analogue of Catlin's conjecture which we now formulate.

\begin{catlinsim}
For any approximating function $\psi\colon \N\to \Rp$

\begin{equation}\label{e:005}
    |\cA_m(\psi)|_m =   1  \quad if \quad
\displaystyle\sum_{r=1}^\infty \ \cN_m(r)  \; \max_{t\ge1} \left(
\dfrac{\psi(rt)}{r t}\right)^m  =\infty \ ,
\end{equation}

\noindent where $\cN_m(r) :=
\#\{(p_1,\dots,p_m)\in\Z^m:(p_1,\dots,p_m,r)=1, \ 0\le p_i<
r\text{ for all }i\}$ is the number of distinct rational points
with denominator $r$ in the unit cube  $\I^m$.
\end{catlinsim}

\noindent In the one dimensional case  $ \cN_m(r) $ is simply
$\varphi(r)$ and (\ref{e:005}) reduces to Catlin's original
conjecture.  For $m \ge 2$, it is not difficult to verify that
 $$
 \cN_m(r) \ \asymp \ r^m \ ,
 $$
where $\asymp$ means comparable and is defined to be the `double' Vinogradov symbol, that
   is both $\ll$ and $\gg$. This together with a nifty geometric argument enables us to
conclude that
$$
\displaystyle\sum_{r=1}^\infty \ \cN_m(r)  \;
\max_{t\ge1} \left( \dfrac{\psi(rt)}{r t}\right)^m
\ \asymp \ \displaystyle \sum_{r=1}^\infty \; \psi(r)^m \ .
$$
The following equivalence is now obvious:

\begin{equation}\label{e:006}
      Corollary \ \ref{corgal}     \hspace*{3ex}  \Longleftrightarrow
  \hspace*{3ex} The \ simultaneous \ Catlin \  conjecture \  for \  m \ge 2.
\end{equation}

The (simultaneous)  higher dimensional analogue of the
Duffin-Schaeffer conjecture \cite{Sprindzuk} requires the stronger
coprimality condition that the coordinates $p_1,\dots,p_m$ of the
vector $\pb\in\Z^m$ are pairwise coprime to $q$. The corresponding
set of $\psi$-approximable points will be denoted by
$\cA_m''(\psi)$ and consists of points $\xb\in \I^m$ for which the
inequality in (\ref{e:004}) is satisfied for infinitely many $(\vv
p,q)\in\Z^m\times\N$ with $(p_j,q)=1$ for $j=1,\dots,m$. In 1990,
Pollington and Vaughan~\cite{PV90} established the {\bf
simultaneous \DS{} conjecture} for $m\ge2$ :

\begin{theorem}[Pollington \&{} Vaughan]
Let $m\ge2$.  For any approximating function $\psi\colon \N\to
\Rp$
$$
|\cA''_{m}(\psi)|_m = \begin{cases} 0
      &\displaystyle\text{if } \ \sum_{r=1}^\infty \varphi(r)^m\left(\frac{\psi(r)}{r}\right)^m<\infty, \\[4ex]
      1
      &\displaystyle\text{if } \  \sum_{r=1}^\infty \varphi(r)^m\left(\frac{\psi(r)}{r}\right)^m=\infty.
                  \end{cases}
$$
\end{theorem}

\noindent Notice that this theorem does not imply Gallagher's
theorem nor does it imply the simultaneous Catlin's conjecture for
$m\ge2$. The books of Sprindzuk~\cite{Sprindzuk} and
Harman~\cite{HarmanMNT} contain a variety of generalisations of the
above results including asymptotic formulae for the number of
solutions, inhomogeneous version of Gallagher's theorem
\cite[Theorem 3.4]{HarmanMNT} and approximation with different
approximating functions in each coordinate.

\section{Dual approximation and Groshev's theorem }
\label{dualgros}

Instead of approximation by rational points as considered
in the previous section, one can consider the closeness of the point $\vv
x=(x_1,\dots,x_n)\in\R^n$ to rational hyperplanes given by the
equations $\vv q\cdot\vv x=p$ with $p\in\Z $ and $\q\in\Z^n$. The point
$\vv x\in\R^n$ will be called \emph{dually $\psi$-approximable}\/ if
the inequality
$$
 |\vv q\cdot \vv x-p|<\psi(|\vv q|)
$$
\noindent holds for infinitely many $(p,\vv q)\in\Z\times\Z^n$,
where $|\vv q|:=|\vv q|_\infty = \max\{|q_1|,\dots,|q_n|\}$. The set
of dually $\psi$-approximable points in $\I^n$ will be denoted by
$\cA^*_n(\psi)$. The argument of the approximating function depends
on the sup norm of $\q$ but it could also easily  be chosen to
depend on $\q\in\Z^n$.  In what follows, we shall consider the sup
norm $|\q|$ and will continue to use $\psi$ to denote an
approximating function on $\N$; i.e. we consider $\psi(|\q|)$ where
$\psi\colon \N\to \Rp$. However in the next section we shall discuss
the general case, when we will use $\Psi$ to denote an approximating
function with argument in $\Z^n$; i.e.  we consider $\Psi(\q)$ where
$\Psi\colon\Z^n\to \Rp$.

The simultaneous and dual forms of approximation are special cases
of a system of linear forms, covered by a  general extension due
to A.~V.~Groshev (see \cite{Sprindzuk}).  This treats real $n
\times m$ matrices $X$, regarded as points in $\R^{nm}$, which are
$\psi$-approximable.  More precisely, $X=(x_{ij}) \in \R^{nm}$ is
said to be $\psi$-approximable if the inequality
$$
  \| \q X \| < \psi(|\q|)
$$
\noindent is satisfied for infinitely many $\q \in \Z^n$. Here $\q
X$ is the system
$$
q_1x_{1j} + \dots + q_n x_{n,j} \hspace*{6ex}  ( 1\leq j \leq m )
$$
of $m$ real linear forms in $n$ variables
and  $\| \q X \| := \max_{1\leq j \leq m } \| \q \cdot X^{(j)}
\|$, where $ X^{(j)}$ is the $j$'th column vector of $X$. As the
set of $\psi$-approximable points is translation invariant under
integer vectors, we can restrict attention to the $nm$-dimensional
unit cube $\I^{nm}$. The set of $\psi$-approximable points in
$\I^{nm}$ will be denoted by
$$
\cA_{n,m}(\psi) := \{X\in \I^{nm}:\|\q X \| < \psi(|\q|) {\text {
for infinitely many }} \q \in \Z^n \}.
$$
Thus, $\cA_m(\psi)=\cA_{1,m}(\psi)$ and
$\cA^*_n(\psi)=\cA_{n,1}(\psi)$. The following result naturally
extends Khintchine's simultaneous theorem to the linear forms
setup.

\begin{theorem}[Groshev]
\label{Groshevthm}  Let $\psi:\N\to\Rp$. Then
$$
    |\cA_{n,m}(\psi)|_{nm} = \begin{cases} 0
      & \ \ds\text{if } \  \sum_{r =1}^\infty r^{n-1}\psi(r)^m<\infty, \\[3ex]
      1
      & \ \ds\text{if }  \  \sum_{r=1}^\infty r^{n-1}\psi(r)^m=\infty \   \text{ and $\psi$ is monotonic}.
                  \end{cases}
$$
\end{theorem}

The counterexample due to Duffin and Schaeffer mentioned in
\S\ref{K} means that the monotonicity condition cannot be dropped
from Groshev's theorem when $m=n=1$. To avoid this situation, let
$m + n > 2$. Then for $n=1$,  we have already seen (Corollary
\ref{corgal}) that as a consequence of Gallagher's theorem the
monotonicity condition can be removed.  Furthermore, the
monotonicity condition can also be removed for $n > 2$ -- this
time due to a  result of Sprindzuk which we discuss in
\S\ref{Sprindzuklf}. However, the $n=2$ situation seems to be
unresolved and we make the following conjecture.

\begin{conjecture}\label{B}
\label{g} Let $\psi:\N\to\Rp$ and suppose that $m + n > 2$. Then
$$
    |\cA_{n,m}(\psi)|_{nm} =
      1
      \qquad\text{if }  \ \  \sum_{r=1}^\infty r^{n-1}\psi(r)^m=\infty \   .
$$
\end{conjecture}

\noindent To reiterate the discussion immediately before the
statement of the conjecture,  it is only the $n=2$ case which is
problematic. It is plausible that it can be resolved using
existing techniques.  Note that for $m +n > 2$, Conjecture \ref{B}
provides a criterion for the size of the `raw' set $
\cA_{n,m}(\psi)$ on which no monotonicity or coprimeness
conditions are imposed.  In view of this, for approximating
functions  with sup norm argument,  Conjecture \ref{B} should
naturally be equivalent to the linear forms analogue of Catlin's
conjecture.

It is possible to formulate the linear forms analogue of both the
Duffin-Schaeffer conjecture and Catlin's conjecture. However, this
will be postponed till the next section in which we consider
multi-variable approximating functions $\Psi\colon\Z^n\to \Rp$ and
thereby formulate the conjectures in full generality. We shall
indeed see that Conjecture \ref{B}  is  equivalent to the linear
forms analogue of Catlin's conjecture for approximating functions
with sup norm argument.

\section{More general approximating functions}\label{Sprindzuklf}

Throughout this section we assume that $n\ge 2$ unless stated
otherwise. In Theorem 13 of Chapter 1 of \cite{Sprindzuk} Sprindzuk
describes a very general setting for a Diophantine system of
linear forms. Let $\{S_{\vv q}\}$ be a sequence of measurable sets
in $\I^{m}$ indexed by integer points $\vv q\in\cZ$, where
$\cZ\subset\Z^n\nz$. Define $\cA_{n,m}(\{S_{\vv q}\})$ to consist
of points $X\in\I^{nm}$ such that there are infinitely many $\vv
q\in\cZ$ satisfying $\vv qX\in S_{\vv q}\pmod 1$. Then

\begin{equation}\label{e:007}
    |\cA_{n,m}(\{S_{\vv q}\})|_{nm}  = \begin{cases} 0 \,
      &\text{if } \   \sum_{\vv q\in\cZ} |S_{\vv q}|<\infty, \\[3ex]
      1 \,
      &\text{if } \   \sum_{\vv q\in\cZ} |S_{\vv q}|=\infty \text{ and any two vectors in $\cZ$ are non-parallel}.
                  \end{cases}
\end{equation}

This extremely general result enables us to generalise  the setup of
\S\ref{dualgros} in two significant ways.  Firstly, we are able to
consider  arbitrary approximating functions $\Psi:\Z^n\to \R^+$
rather than restrict the argument of $\Psi$ to the sup norm of $\vv
q$. Secondly, we are naturally able to consider inhomogeneous
problems. Define

\begin{equation}\label{e:008}
  \cA_{n,m}^{\vv b}(\Psi) := \{X\in \I^{nm}:\|\q X +\vv
b \| < \Psi(\q) {\text { for
    infinitely many }} \q \in \Z^n \}\,.
  \end{equation}

\noindent Here $\vv b\in[0,1)^m$ is a fixed vector that represents the
`inhomogeneous' or `shifted' part of approximation. Now let $\vv
P^n$ denote the set of \emph{primitive}\/  vectors in $\Z^n$ --
non-zero integer vectors with coprime components. It is easy to
see that the statement given by  (\ref{e:007}) with $\cZ = \vv
P^n$ specializes to give the following result -- essentially
Theorem 14 in~\cite{Sprindzuk}.

\begin{theorem}[Sprindzuk]
\label{Sprind} Let $\Psi\colon\Z^n\to \Rp$ and suppose  that
$\Psi(\q) = 0 $ for $\q \notin \vv P^n$.  Then for $n \ge 2$
  \begin{equation}\label{e:009}
    |\cA^{\vv b}_{n,m}(\Psi)|_{nm} = \begin{cases}  \ 0 \
      &\text{if } \ \sum_{\vv q\in \Z^n} \Psi(\vv q)^m<\infty \ , \\[2ex]
       \ 1 \
      &\text{if }  \ \sum_{\vv q\in\Z^n} \Psi(\vv q)^m=\infty \
      .
                  \end{cases}
  \end{equation}
\end{theorem}

\noindent Of course the primitivity condition  in the theorem,
namely that $\Psi$ vanishes on non-primitive integer vectors,
imposes a primitivity condition on the set $\cA^{\vv b}_{n,m}(\Psi)$
-- namely that the vectors $\q \in \Z^n$ associated with
(\ref{e:008}) are primitive.

We now consider a special case of Sprindzuk's theorem in which the
argument of $\Psi$ is restricted to the sup norm. In keeping with
the notation used in \S\ref{dualgros}, we write $\psi$ for $\Psi$
and so  $\Psi(\vv q)=\psi(|\vv q|)$ for $\vv q$ in $\Z^n$.  Let
$n\ge 3$. Then the number of primitive vectors $\vv q$ in $\Z^n$
with $|\vv q|=r$ is comparable to $r^{n-1}$. Thus the number of
primitive vectors is comparable to the number of vectors without
any primitivity restriction.  It follows that the divergence
condition in (\ref{e:009}) is equivalent to $ \sum_{r=1}^\infty
r^{n-1}\psi(r)^m=\infty $. The latter is precisely the divergent
sum appearing in Groshev's theorem.  The upshot of this is that
Sprindzuk's theorem  removes the monotonicity requirement
from Groshev's theorem when $n \ge 3$. A similar argument in the
case $n=2$ does not yield an equivalent improvement in Groshev's
theorem. The reason for this is simply the fact that the number of primitive
vectors $\vv q$ in $\Z^2$ with $|\vv q|=r$ is comparable to
$\varphi(r)$ whereas the number of vectors without any primitivity
restriction is comparable to $r$. In short, when $n=2$ the
divergent sum appearing in Sprindzuk's theorem is not equivalent
to that appearing in Groshev's theorem.

The primitivity condition cannot in general be omitted from
Sprindzuk's theorem. To give a counterexample, we  consider the
case that  $m=1$ and $n \ge 2$. For $\vv
q=(q_1,\dots,q_n)\in\Z^n\nz$ let
\begin{equation}\label{e:010}
    \Psi_\vartheta(\vv q) := \left\{\begin{array}{cl}
      \vartheta(|q_1|) & \text{if }\vv q=(q_1,0,\dots,0) , \\[2ex]
      0 & \text{otherwise,}
    \end{array}\right.
\end{equation}
where $\vartheta$ is the function constructed by Duffin
and Schaeffer (see \S\ref{K}). Obviously $$\sum_{\vv q\in\Z^n\nz}
\Psi_\vartheta(\vv q)=\sum_{q\in\Z\nz} \vartheta(|q|)=\infty.$$ On the other
hand, $\cA_{n,1}(\Psi)=\cA(\vartheta)\times\I^{n-1}$ so that
$|\cA_{n,1}(\Psi_\vartheta)|=|\cA(\vartheta)|\cdot|\I^{n-1}|=0\cdot1=0$.
Clearly $\Psi_\vartheta$ does not satisfy the primitivity condition in Sprindzuk's theorem.

The above counterexample implies that the primitivity condition cannot be omitted from Sprindzuk's
theorem for the dual form of approximation; namely when considering the dual set  $\cA_{n}^*(\Psi) := \cA_{n,1}(\Psi)$.
 However, if $m > 1$ so that  we are dealing with a system of more than one linear
 form, no similar counterexample appears to be possible.
 Indeed we strongly believe in the truth of the following conjecture concerning the set
 \begin{equation}\label{e:011}
 \begin{array}[b]{rl}
 \cA'_{n,m}(\Psi) := \{X\in \I^{nm}:& |\q X+\vv p| <
                 \Psi(\q) {\text { for infinitely many }} (\vv p,\vv q) \in
 \Z^m\times\Z^n\\[1ex] &\text{ with }(p_1,\dots,p_m,q_1,\dots,q_n)=1 \}\end{array}
   \end{equation}

\begin{conjecture}\label{C}
Let $\Psi:\Z^n\to\Rp$ and suppose that $m>1$. Then
 $$
     |\cA'_{n,m}(\Psi)|_{nm} =
       1 \qquad
       \text{if } \ \   \sum_{\vv q\in\Z^n\nz} \Psi(\vv q)^m=\infty \ .
 $$
\end{conjecture}

Obviously for every $\vv q\in\vv P^n$ the coprimality condition in
(\ref{e:011}) is satisfied. Therefore
$\cA_{n,m}(\Psi)=\cA'_{n,m}(\Psi)$ for any $\Psi$ vanishing outside
of $\mathbf{P}^n$. Thus Conjecture~\ref{C} covers Sprindzuk's
theorem in the case of $m>1$. Furthermore, since
$\cA'_{n,m}(\Psi)\subset\cA_{n,m}(\Psi)$, Conjecture~\ref{C} would
imply the following statement for the `raw' set $ \cA_{n,m}(\Psi)$
on which no monotonicity or coprimeness conditions are imposed.
\begin{conjecture}\label{D}
Let $\Psi:\Z^n\to\Rp$ and suppose that $m>1$. Then
$$
     |\cA_{n,m}(\Psi)|_{nm} =
       1 \qquad
       \text{if } \ \   \sum_{\vv q\in\Z^n\nz} \Psi(\vv q)^m=\infty \ .
$$
\end{conjecture}

\noindent As one should expect, we will see below that Conjecture \ref{D} is naturally  equivalent to the linear forms
analogue of Catlin's conjecture for $m > 1$. Recall that the above conjectures, are straightforwardly
established in the complementary convergent cases.

The counterexample given by (\ref{e:010}) shows that
Conjectures~\ref{C} and \ref{D} are not valid when $m=1$; i.e.  the dual form of approximation.
We now deal with the case $m=1$. It is relatively easy to show that the set
$\cA'_{n,1}(\Psi)$ has measure zero if
the sum
\begin{equation}\label{e:012}
    \sum_{\vv q\in\Z^n\nz}\frac{\varphi\big(\gcd(\vv q)\big)}{\gcd(\vv q)}
\ \Psi(\vv q)\ = \ \sum_{d=1}^\infty\frac{\varphi(d)}{d} \sum_{\vv q'\in\,\vv P^n} \Psi(d\vv q')
\end{equation}
converges. Here $\gcd(\vv q)$ denotes the greatest common divisor of the
components of $\vv q\in\Z^n$. The following can
be regarded as a generalisation of the Duffin-Schaeffer conjecture
to the case of dual approximation.
\begin{DSdual}
Let $\Psi:\Z^n\to\R^+$ and suppose that $m=1$. Then
$$
     |\cA'_{n,1}(\Psi)|_{n} =
       1 \qquad
       \text{if \ \ \ $(\ref{e:012})$ \ diverges} .
$$
\end{DSdual}
It is clear that this conjecture includes the original
Duffin-Schaeffer conjecture. It would be desirable to find natural
conditions on $\Psi$ which make the conjecture genuinely multi-dimensional.
For example, in the genuine multi-dimensional case
it is natural to exclude  approximating functions $\Psi$ like  $\Psi_\vartheta$ given by  (\ref{e:010}).
The hope is that the genuine multi-dimensional case is
`easier' than the one dimensional case.
Recall that in the case of simultaneous approximation, the multi-dimensional
Duffin-Schaeffer conjecture has been proved.
Also notice that if $n\ge2$ and there exists  $d\in\N$ such that
$$
\sum_{\vv q'\in\,\vv P^n} \Psi(d\vv q')\,,
$$
the internal sum in (\ref{e:012}), diverges,  then
the conjecture is reduced to Sprindzuk's theorem. Therefore, it is also reasonable to assume
that the internal sum in (\ref{e:012}) is  finite irrespective of $d$.

Regarding the `raw' set $\cA_{n,1}(\Psi)$ on which no monotonicity
or coprimeness conditions are imposed, it is natural to formulate
the analogue of Catlin's conjecture. For $\vv q\in\Z^n$ let
$$
\cN^*_n(\vv q) \ := \ \#\{0< p\le|\vv q|:(p,\vv
q)=1\} .
$$
It is easy to see that
$$
\cN^*_n(\vv q) \ \asymp \ \frac{\varphi\big(\gcd(\vv q)\big)}{\gcd(\vv q)}\ |\vv
q|
$$
and moreover that the set $\cA_{n,1}(\Psi)$ has measure zero if
the sum
\begin{equation}\label{e:013}
    \sum_{\vv q\in\Z^n\nz}\cN^*_n(\vv q)\ \max_{t\ge 1}
\frac{\Psi(t\vv q)}{t|\vv q|}
\end{equation}
converges. The following can
be regarded as a generalisation of Catlin's conjecture
to the case of dual approximation:
\newtheorem{Cdual}{The dual Catlin conjecture : }
\renewcommand{\theCdual}{}
\begin{Cdual} Let $\Psi:\Z^n\to  \R^+$ and suppose that $m=1$.
Then
$$
     |\cA_{n,1}(\Psi)|_{n} =
       1 \qquad
       \text{if \ \ \ $(\ref{e:013})$ \ diverges} .
$$
\end{Cdual}

\noindent For $n\ge2$ and $\Psi(\vv q)=\psi(|\vv q|)$ the above conjecture is equivalent
to Conjecture~\ref{B} with $m=1$. Recall, that for $n\ge3$ the latter is known to be true.

\medskip

Finally, for the sake of completeness,
we extend the above dual conjectures to the general linear forms setup.
For the analogue of the
Duffin-Schaeffer conjecture it is natural to impose a
coprimality condition on each linear form. Let
$$
\cA''_{n,m}(\Psi) := \{X\in \I^{nm}:\begin{array}[t]{r}
                |\q X+\vv p| <
                \Psi(\q) {\text { for infinitely many }} (\vv p,\vv q) \in
\Z^m\times\Z^n\\[1ex] \text{ with }(p_j,q_1,\dots,q_n)=1\text{ for every }j=1,\dots,m
\}\,.\end{array}
$$
  It is easy to show that the set
$\cA''_{n,m}(\Psi)$ has measure zero if
the sum
  \begin{equation}\label{e:014}
   \sum_{\vv q\in\Z^n\nz}\left(\frac{\varphi\big(\gcd(\vv q)\big)}{\gcd(\vv
q)}\,\Psi(\vv q)\right)^m
  \end{equation}
  converges. This leads to the following complementary problem.

\begin{DSsystems} Let $\Psi:\Z^n\to  \R^+$.
Then
  $$
    |\cA''_{n,m}(\Psi)|_{nm} = 1
      \qquad \text{if the sum $(\ref{e:014})$ \ diverges.}
  $$
\end{DSsystems}

Regarding the `raw' set $\cA_{n,m}(\Psi)$ on which no monotonicity
or coprimeness conditions are imposed, we formulate
the analogue of Catlin's conjecture.
For $\vv q\in\Z^n$ let
\begin{equation}\label{e:015}
    \cN_{n,m}(\vv q) \ := \ \#\{\vv p\in\Z^m: 0\le |\vv p|\le|\vv q|,\ (\vv p,\vv
q)=1\} .
\end{equation}
It can be verified that the set $\cA_{n,m}(\Psi)$ has measure zero
if the sum
\begin{equation}\label{e:016}
    \sum_{\vv q\in\Z^n\nz}\cN_{n,m}(\vv q)\ \left(\max_{t\ge 1}
\frac{\Psi(t\vv q)}{t|\vv q|}\right)^m
\end{equation}
converges. This leads to the following complementary problem.
\begin{catlinlin} Let $\Psi:\Z^n\to  \R^+$.
Then
$$
     |\cA_{n,m}(\Psi)|_{nm} =
       1 \qquad
       \text{if \ \ \ $(\ref{e:016})$ \ diverges} .
  $$
\end{catlinlin}
On modifying the argument that enables us to establish the
equivalence (\ref{e:006}) within the simultaneous setup, we
obtain the following statements in which the divergence of
(\ref{e:016}) is simplified.

{\it
$$
~ \hspace*{-8ex}
\text{Conjecture~\ref{D}}\qquad\Longleftrightarrow\qquad\text{Linear
forms Catlin's conjecture for $m\ge2$}\,.
$$
}

$$
 \begin{array}[b]{rl}
Conjecture~\ref{B} \qquad\Longleftrightarrow\qquad  Linear \ forms&
\!\!\! Catlin's \ conjecture \ for \  m+n>2  \\ & and  \ \Psi(\vv
q)=\psi(|\vv q|) \, .  \end{array}
$$

\section{The theorems of \Ja{} and Besicovitch}\label{J}

The results of \S\ref{K}--\ref{Sprindzuklf} can be regarded as the
\emph{probabilistic theory of Diophantine approximation}. Indeed,
these results indicate the probability of a certain Diophantine
property and include both qualitative results like Khintchine's
theorem and their quantitative versions (see
\cite{HarmanMNT,Sprindzuk}). Furthermore, the results are rigid in
the sense that the indicated probability is always either {\it
zero} or {\it one}. Even in the case of the profound and as yet
unsolved problem of Duffin-Schaeffer, it is known (see
\cite{HarmanMNT,Sprindzuk}) that the measure of $\cA'(\psi)$ and
indeed $\cA(\psi)$ must satisfy this rigid `zero-one' law.

As the results considered obey zero-one laws, they always involve
`exceptional' sets of measure zero. The probabilistic theory of
Diophantine approximation doesn't tell us anything more about the
`size' of these exceptional sets, although it is intuitively clear
that it should depend on the choice of the approximating function.
This leads us to a more delicate study which makes use of various
concepts from geometric measure theory -- in particular Hausdorff
measure and dimension.

\subsection{Hausdorff measures and dimension}

In what follows, a {\em dimension function\,} $f:\Rp \to \Rp $ is a
left continuous, monotonic function such that $f(0)=0$.
Suppose $F$ is a  subset  of $\R^n$. Given a ball $B$ in $\R^n$, let
$r(B)$ denote the radius of $B$. For $\rho
> 0$, a countable collection $ \left\{B_{i} \right\} $ of balls in
$\R^n$ with $r(B_i) \leq \rho $ for each $i$ such that $F \subset
\bigcup_{i} B_{i} $ is called a {\em $ \rho $-cover for $F$}. Define
$$
 {\cal H}^{f}_{\rho} (F) \, := \, \inf \ \sum_{i} f(r(B_i)),
$$
where the infimum is taken over all $\rho$-covers of $F$. The {\it
Hausdorff $f$--measure of $F$} denoted by $ {\cal H}^{f} (F)$ is
defined as
$$ {\cal H}^{f} (F) := \lim_{ \rho \rightarrow 0} {\cal
H}^{f}_{\rho} (F) \; = \; \sup_{\rho > 0 } {\cal H}^{f}_{\rho} (F)
\; .
$$
In the case that  $f(r) = r^s$ ($s \ge 0$), the measure $ \cH^f $
is the more common {\em$s$-dimensional Hausdorff measure} $\cH^s
$, the measure $\cH^0$ being the cardinality of $F$.  Note that
when $s$ is a positive integer, $\cH^s$ is a constant multiple of
Lebesgue measure in $\R^s$ and that when $s=1$, the measures
coincide.  Thus if the $s$-dimensional Hausdorff measure of a set
is known for each $s>0$, then so is its $n$-dimensional Lebesgue
measure for each $n\ge1$. The following easy property
$$
\cH^{s}(F)<\infty\quad  \Longrightarrow\quad
\cH^{s'}(F)=0\qquad\text{if }s'>s
$$
implies that there is a unique real point $s$ at which the Hausdorff
$s$-measure drops from infinity to zero (unless the set $F$ is
finite so that $\cH^s(F)$ is never infinite). This point is called
the Hausdorff dimension of $F$ and is formally defined as
$$
\dim F :=  \inf \left\{ s>0 : \cH^{s} (F) =0 \right\}.
$$

\noindent The Hausdorff dimension has been established  for many
number theoretic sets, e.g.  $\cA(\tau)$ (this is the \JB{}
theorem discussed below), and is easier than determining
the Hausdorff measure. Further details regarding Hausdorff measure
and dimension can be found in~\cite{FalcGFS,MattilaGS}.

\vspace{3ex}

\subsection{The theorems}

The first step towards the study of Hausdorff measure of the set
of $\psi$-approximable points was made by Jarn\'{\i}k ~\cite{Ja29}
in 1929 and independently by Besicovitch ~\cite{Bes34} in 1934.
They determined the Hausdorff dimension of the set $\cA(q\mapsto
q^{-\tau})$, usually denoted by $\cA(\tau)$, where $\tau>0$.
\begin{theorem}[\Ja-Besicovitch]
$$
  \dim \cA(\tau) = \begin{cases}
       \ \frac2{\tau + 1} \; & \ds\text{if } \ \
      \tau>1 \; ,\\[2ex]
       \ 1 \: & \ds\text{if } \ \ \tau \leq 1 \; .
    \end{cases}
$$
\end{theorem}

Note that for $\tau \leq 1$ the result is trivial since $\cA(\tau) =
\I$ as a consequence of Dirichlet's theorem.  Thus the main content
is when $\tau > 1$. In this case, the dimension result implies that
$$ \cH^s \left(\cA(\tau)\right)=\left\{\begin{array}{ll} 0 & \ \ \ \
{\rm if} \;\;\; s \; > \; \frac2{\tau + 1}\,,   \\ [2ex] \infty & \
\ \ \ {\rm if} \;\;\; s \; < \; \frac2{\tau + 1}\,,
\end{array}\right. $$
but gives no information regarding the $s$--dimensional Hausdorff
measure of $\cA(\tau)$ at the critical value $s=\dim \cA(\tau)$.
In a deeper study, \Jarnik~\cite{Ja31} essentially established the
following general Hausdorff {\em measure} result for simultaneous
\DA.

\begin{theorem}[Jarn\'{\i}k]
\label{jarmeas} Let $f$ be a dimension function such that $r^{-m}
\, f(r)\to \infty$ as $r\to 0 \, $ and $r^{-m} \, f(r) $ is
decreasing. Then
$$
\cH^f\left(\cA_m(\psi)\right)=\left\{\begin{array}{cl} 0 &
\ds\text{if } \;\;\; \sum_{r=1}^{\infty}  \ \; r^m \,
f\left(\frac{\psi(r)}r\right)^m
 <\infty \; ,\\[3ex]
\infty & \ds\text{if } \;\;\; \sum_{r=1}^{\infty} \  \; r^m \,
f\left(\frac{\psi(r)}r\right)^m =\infty   \  \text{ and $\psi$ is
                    monotonic}.
\end{array}\right.$$
\end{theorem}

 With $m=1$ and $f(r)=r^s$, Jarn\'{\i}k's theorem not only
gives the above dimension result but implies that
\begin{equation}\label{e:017}
 \cH^{\frac2{\tau + 1}} \left(\cA(\tau)\right) =
\infty \quad \ds{\rm \ if} \ \
      \tau>1.
\end{equation}
\noindent For monotonic approximating functions $\psi:\N\to\Rp$,
Jarn\'{\i}k's theorem provides a beautiful and simple criteria for
the `size' of the set $\cA_m(\psi)$ expressed in terms of
Hausdorff measures. Naturally, it can be regarded as the Hausdorff
measure version of Khintchine's simultaneous theorem. As with the
latter, the divergence part constitutes the main substance.
Notice, that the case when $ \cH^f $ is comparable to
$m$--dimensional Lebesgue measure  (i.e. $ f(r)= r^m$) is excluded
by the condition $r^{-m} \, f(r)\to \infty$ as $r\to 0 \, $.
Analogous to Khintchine's original statement, {\em in
Jarn\'{\i}k's original statement the additional hypotheses that $r
\p(r)^m$ is decreasing, $r \p(r)^m \to 0 $ as $ r \to \infty $ and
that $r^{m+1} f(\p(r)/r) $ is decreasing were assumed.} Thus, even
in the simple case when $m=1$, $f(r) = r^s $ $ (s\ge 0) $  and the
approximating function is given by $\psi(r) = r^{-\tau} \log r $
$(\tau > 1)$, Jarn\'{\i}k's original statement gives no
information regarding the $s$--dimensional Hausdorff measure  of
$\cA(\psi)$ at the critical exponent $s=2/(\tau+1)=\dim\cA(\psi)$.
That this is the case is due to the fact that $r^2f(\p(r)/r) $ is
not decreasing. Recently, however, it has been shown
in~\cite{BDV06} that the monotonicity of $\psi$ suffices in
\Jarnik's theorem. In other words, the additional hypotheses
imposed by Jarn\'{\i}k are unnecessary. Furthermore, with the
theorems of Khintchine and Jarn\'{\i}k as stated above it is
possible to combine them to obtain a single unifying statement.

\begin{theorem}[Khintchine-Jarn\'{\i}k]
\label{khijarmeas} Let $f$ be a dimension function such that
 $r^{-m} \, f(r) $
is monotonic. Then
$$
\cH^f\left(\cA_m(\psi)\right)=\left\{\begin{array}{cl} 0 &
\ds\text{if } \;\;\; \sum_{r=1}^{\infty}  \ \; r^m \,
f\left(\frac{\psi(r)}r\right)^m
 <\infty \; ,\\[3ex]
\cH^f(\I^m) & \ds\text{if } \;\;\; \sum_{r=1}^{\infty} \  \; r^m
\, f\left(\frac{\psi(r)}r\right)^m =\infty   \  \text{ and $\psi$
is monotonic}.
\end{array}\right.$$
\end{theorem}

For monotonic approximating functions, the Khintchine-Jarn\'{\i}k
theorem provides a complete measure theoretic description of
$\cA_m(\psi)$. Clearly, when $f(r) = r^m$ the theorem corresponds
to Khintchine's theorem. It would be quite natural to suspect that
such a unifying statement is established by combining two
independent results: the Lebesgue measure statement (Khintchine's
theorem) and the Hausdorff measure statement (Jarn\'{\i}k's
theorem). Indeed, the underlying method of proof of the individual
statements are dramatically different.  However, this is not the
case. In view of the Mass Transference Principle recently
established in~\cite{BV06} one actually has that
\begin{center} Khintchine's Theorem $ \hspace{4mm} \Longrightarrow
\hspace{4mm} $ Jarn\'{\i}k's Theorem.  \end{center} Thus, the
Lebesgue theory of $\cA_m(\psi)$ underpins the general Hausdorff
theory. At first glance this is rather surprising because the
Hausdorff theory had previously  been  thought to be a
subtle refinement of the Lebesgue theory. Nevertheless, the Mass
Transference Principle allows us to transfer Lebesgue measure
theoretic statements for limsup sets to Hausdorff statements and naturally obtain
a complete metric theory. That this is the case is by no means a
coincidence -- see \cite{BV06,BV06Slicing, BVparis} for various
points of view.

\section{Mass Transference Principle}\label{mtp0}

Given a dimension function $f$, define the following transformation
on balls in $\R^m$:
$$
\textstyle B=B(x,r)\mapsto B^f:=B(x,f(r)^{1/m}) \ .
$$
When $f(x)=x^s$ for some $s>0$ we also adopt the notation $B^s$ for
$B^f$. Clearly $B^m=B$. Recall that $\cH^m$ is comparable to the
$m$-dimensional Lebesgue measure. The limsup of a sequence of
balls $B_i$, $i=1,2,3,\ldots$ is
$$
\limsup_{i\to\infty}B_i:=\bigcap_{j=1}^\infty\ \bigcup_{i\ge j}B_i \
.
$$

\noindent For such limsup sets, the following statement (the \MTP{})
is the key to obtaining Hausdorff measure statements from Lebesgue
statements.

\begin{theorem}[Beresnevich \&{} Velani]
Let $\{B_i\}_{i\in\N}$ be a sequence of balls in $\R^m$ with
$\diam(B_i)\to 0$ as $i\to\infty$. Let $f$ be a dimension function
such that $x^{-m}f(x)$ is monotonic.  For  any finite  ball
$B$ in $\R^m$, if
$$
\cH^m\big(\/B\cap\limsup_{i\to\infty}B^f_i{}\,\big)=\cH^m(B) \
$$
then
$$
\cH^f\big(\/B\cap\limsup_{i\to\infty}B^m_i\,\big)=\cH^f(B) \ .
$$
\end{theorem}

There is one point that is well worth making. The Mass
Transference Principle is purely a statement concerning limsup
sets arising from a sequence of balls. There is absolutely no
monotonicity assumption on the radii of the balls.
Even the imposed condition that $\diam(B_i)\to0$ as
$i\to\infty$ is redundant but is included to avoid unnecessary
tedious discussion.

\subsection{A Hausdorff measure Duffin-Schaeffer conjecture}
\label{subsec:HMDSC}

As an application of the Mass Transference Principle, we shall see
that the simultaneous Duffin-Schaeffer conjecture implies the
corresponding conjecture for Hausdorff measures.

Let $f$ be a dimension function.
 A straightforward covering argument
making use of the $\limsup$ nature of $\cA''_m(\psi)$ implies that
$$
 \cH^f(\cA''_m(\psi)) = 0  \ \ \ \ \ {\rm if \ } \ \ \ \ \
 \sum_{q=1}^{\infty}  f\left(\frac{\psi(q)}{q}\right) \; \varphi(q)^m
\ < \ \infty .
$$
It is therefore natural to make the following conjecture (see \cite{BV06}) which
can be regarded as the simultaneous Duffin-Schaeffer conjecture for Hausdorff measures.
\begin{conjecture}\label{DSHM}
Let $f$ be a dimension function such that $r^{-m} f(r)$ is monotonic. Then
$$
\cH^f(\cA''_m(\psi)) = \cH^f(\I^m)\qquad \text{if }
\qquad\sum_{q=1}^{\infty} f\left(\frac{\psi(q)}{q}\right) \;
\varphi(q)^m  = \infty \,.
$$
\end{conjecture}

When $\psi$ is monotonic, Conjecture~\ref{DSHM} reduces to the Khintchine-\Jarnik{}
theorem. It turns out that Conjecture~\ref{DSHM}, a
refinement of the Duffin-Schaeffer problem, is simply its
consequence~\cite{BV06}.

\begin{theorem}[Beresnevich \& Velani]\label{onemore}
$$
\text{The simultaneous Duffin-Schaeffer conjecture\qquad $\Longleftrightarrow$ \qquad Conjecture~\ref{DSHM}}
$$
\end{theorem}

Conjecture~\ref{DSHM} contains the simultaneous Duffin-Schaeffer conjecture.
In order to prove the converse note that
$\cA''_m(\psi)$ is the limsup set of the sequence of balls given by
$$
\text{$|qx-\vv p|<\psi(q)$ with $(q,\vv p)\in\N\times\Z^m$ and $0\le p_j\le q$ for all $j=\overline{1,m}$.}
$$
First we can dispose of the case that $\psi(q)/q \nrightarrow 0 $,
as $q \to \infty$ as otherwise the result is trivial. We are given
that $ \sum f( \psi(q)/q ) \; \varphi(q)^m = \infty $. Let $
\theta(q) := q \, f( \psi(q)/q )^{1/m}$. Then $\theta$ is an
approximating function and $ \sum (\varphi(q) \, \theta(q)/q )^m =
\infty $. Thus, on using the supremum norm, the Duffin-Schaeffer
conjecture implies that $ \cH^m(B \cap \cA''_m(\theta)) =
\cH^m(B)$ for any ball $B$ in $\I^m$. It now follows via the Mass
Transference Principle with $B = \I^m$ that $ \cH^f(\cA''_m(\psi)) = \cH^f(\I^m)$
and establishes Theorem~\ref{onemore}. Since the simultaneous
Duffin-Schaeffer conjecture is know to be true for $m\ge 2$,
Theorem~\ref{onemore} implies the following result.

\begin{corollary}
Conjecture~\ref{DSHM} holds for $m\ge2$.
\end{corollary}

\medskip

In a similar fashion, the Mass Transference Principle yields the following generalisation of
Gallagher's theorem.

\begin{theorem} \label{galyest}
Let $m\ge2$. Let $f$ be a dimension function such that $r^{-m} f(r)$ is
monotonic. Then
$$
    \cH^f(\cA'_{m}(\psi)) = \begin{cases} 0 \
      &\ds\text{if }  \quad \sum_{r=1}^\infty f\left(\frac{\psi(r)}{r}\right) \; r^m<\infty \, , \\[3ex]
      \cH^f(\I^m) \
      &\ds\text{if }   \quad \sum_{r=1}^\infty f\left(\frac{\psi(r)}{r}\right) \;
      r^m=\infty \,
      .
                  \end{cases}
$$
\end{theorem}

\noindent Note that since $\cA_m'(\psi)\subseteq \cA_m(\psi)$, Theorem \ref{galyest}
implies the divergent part of the Khintchine-\Jarnik{} theorem. Furthermore,
we deduce that \textit{the monotonicity
condition in the Khintchine-\Jarnik{} theorem is redundant if $m\ge2$}.

It is remarkable that by using the \MTP{} one can deduce the
\Jarnik-Besicovitch theorem from Dirichlet's theorem. Moreover,
one obtains the stronger measure statement given by
(\ref{e:017}). Finally, we point out that all the Lebesgue
measure statements on simultaneous Diophantine approximation can
be generalised to the Hausdorff measure setting as above.

\section{Mass Transference Principle for systems of linear
forms}\label{MTPSLF}

The Mass Transference Principle of \S\ref{mtp0} deals with
$\limsup$ sets which are defined as a sequence of balls. However,
the `slicing' technique introduced in \cite{BV06Slicing}  extends
the \MTP{} to deal with $\limsup$ sets defined as a sequence of
neighborhoods of `approximating' planes.  This naturally enables
us to generalise the Lebesgue measure statements of
\S\ref{dualgros}-\ref{Sprindzuklf} for systems of linear forms to
Hausdorff measure statements. In particular,  Groshev's theorem
can be extended to obtain a linear forms  analogue of the
Khintchine--\Jarnik{} theorem.

Throughout  $k, m \ge 1 $ and $ l\ge0$ are integers such that
$k=m+l$. Let $\cR=(\ra )_{\alpha \in J}$ be a family of planes in
$\R^k$ of common dimension $l$ indexed by an infinite countable
set $J$. For every $\alpha\in J$ and $\delta\ge 0$ define
$$ \Delta(R_\alpha,\delta) := \{\vv x \in \R^k: \dist(\vv
x,R_\alpha) < \delta\} \ . $$ Thus $\Delta(R_\alpha,\delta)$ is
simply the $\delta$--neighborhood of the $l$--dimensional plane
$R_\alpha$. Note that by definition, $ \Delta(R_\alpha,\delta) =
\emptyset $ if $\delta =0$. Next, let
$$\Upsilon : J \to \Rp : \alpha\mapsto
\Upsilon(\alpha):=\Upsilon_\alpha$$ be a non-negative, real valued
function on $J$. We assume that for every $\epsilon
>0$ the set $\{\alpha\in J:\Upsilon_\alpha>\epsilon \}$ is finite.
This condition implies that $\Upsilon_\alpha  \to 0 $ as $\alpha$
runs through $J$. We now consider the following `$\limsup$' set,
$$ \La(\Upsilon):=\{\vv x\in\R^k:\vv
x\in\Delta(R_\alpha,\Upsilon_\alpha)\ \mbox{for\ infinitely\ many\
}\alpha\in J\} \ . $$ Note that in view of the conditions imposed on
$k,l$ and $m$ we have that $l < k$. Thus the dimension of the
`approximating' planes $ R_\alpha$ is strictly less than that of the
ambient space $\R^k$. The situation when $l =k $ is of little
interest.

The following statement is a generalisation of the \MTP{} to the
case of systems of linear forms.

\begin{theorem}[Beresnevich \&{} Velani]\label{MTPLF}
Let $\cR$ and $\Upsilon$ as above be given. Let $V$ be a linear
subspace of\/ $\R^k$ such that $\dim V=m=\mathrm{codim}\,\cR$ and

~ \qquad  \qquad $(i)$\quad \ $V \ \cap \  R_\alpha \ \neq \
\emptyset $ \quad for all $ \ \alpha\in J \ $,

~ \qquad  \qquad $(ii)$\quad $\sup_{\alpha\in J}\diam( \,
V\cap\De(R_\alpha,1) \, ) \ < \ \infty \ $ .

\noindent Let $f$ and $g : r \to g(r):= r^{-l} \, f(r)$ be
dimension functions such that $r^{-k}f(r)$ is monotonic.
For any finite ball $B$ in $\R^k$, if
\begin{equation}\label{e:018}
    \cH^k \big(
\, B  \cap \La \big(g(\Upsilon)^{\frac1m} \big) \,  \big) \, = \,
\cH^k(B)
\end{equation}
 then
\begin{equation}\label{e:019}
    \cH^f \big( \, B \cap \Lambda(\Upsilon) \, \big) \, = \, \cH^f(B) \, .
\end{equation}
\end{theorem}

\medskip

When $l=0$,  so that  $\cR$ is a collection of points in $\R^k $,
conditions (i) and (ii) are trivially satisfied. When $l \ge 1 $,
so that $\cR$ is a collection of $l$--dimensional planes  in $\R^k
$, condition (i) excludes planes $\ra$  parallel to $V$ and
condition (ii) simply means that the angle at which $R_\alpha$
`hits'  $V$ is bounded away from zero by a fixed constant
independent of $\alpha \in J$. This in turn implies that each
plane in $\cR$ intersects $V$ at exactly one point. The upshot is
that the conditions  (i) and (ii) are not particularly
restrictive. However, we believe that they are actually redundant.

\begin{conjecture}
Theorem~\ref{MTPLF} is valid without hypothesis (i) and (ii).
\end{conjecture}

\bigskip

As an application of the  \MTP{} for systems of linear forms, we
shall obtain the following Hausdorff measure generalisation of
Sprindzuk's theorem.

\begin{theorem}\label{bv100}
 Let $\Psi:\Z^n\to\Rp$ and
suppose  that $\Psi(\vv q) = 0 $ for $\vv  q \notin \vv P^n$. Let
$f$ and $g : r \to g(r):= r^{-m(n-1)} \, f(r)$ be dimension
functions such that $r^{-mn}f(r)$ is monotonic. Then for $n \ge 2$
$$
\cH^f(\cA_{n,m}^{\vv b}(\Psi)) = \begin{cases} 0 \
&\ds\text{if } \qquad \sum_{\vv q\in\Z^n\nz}g\left(\frac{\Psi(\vv q)}{|\vv q|}\right)\ |\vv q|^m<\infty, \\[4ex]
\cH^f(\I^{nm}) \ &\ds\text{if } \qquad \sum_{\vv
q\in\Z^n\nz}g\left(\frac{\Psi(\vv q)}{|\vv q|}\right) \ |\vv
q|^m=\infty .
                  \end{cases}
$$
\end{theorem}

The convergence case is readily established using standard
covering arguments. We will concentrate on the divergence case and assume
that $r^{-k}f(r)$ is decreasing. The statement is almost obvious if the latter
is not the case. When the sum given in the theorem
diverges, there is a $j\in\{1,\dots,n\}$ such that $$\sum_{\vv
  q\in\Z^n\nz}g\left(\frac{\Psi_j(\vv q)}{|\vv q|}\right) \, |\vv
q|^m \ = \ \infty \, , $$  where $\Psi_j(\vv q)$ vanishes on
$\vv q$ with $|\vv q|\not=|q_j|$ and equals $\Psi(\vv q)$
otherwise. Fix such a $j$. For each point $X\in\cA_{n,m}^{\vv
b}(\Psi_j)$ there are infinitely many $\vv q\in\Z^n\nz$ such
that
\begin{equation}\label{e:020}
    \|\vv qX+\vv b\|<\Psi_j(\vv q)\,.
\end{equation}
In fact, we have $|\vv q|=|q_j|$ for every solution $\vv q$ of
(\ref{e:020}). Now, let
$$
J := \{ (\vv q, \vv p) \in \Z^n\setminus\{\vv 0\}  \times \Z^m  :
|\vv q | = |q_j| \} ,
$$
$\alpha := (\vv q, \vv p) \in J $, $R_\alpha
:= R_{{\vv q},{\vv p}}$ where
$$
 R_{\vv q,\vv p}:=\{X\in\R^{nm}:\vv qX+\vv p+\vv b=0\}\,,
$$
and $ \Upsilon_\alpha := \textstyle{\frac{\Psi_j(\vv
q)}{\sqrt{{\vv q}.{\vv q}}} } $. Then,
\begin{equation}\label{e:021}
    \cA_{n,m}^{\vv b}(\Psi)  \  \supset \ \cA_{n,m}^{\vv
b}(\Psi_j) \ = \ \La(\Upsilon)  \, \cap \, \I^{nm} \ \ .
\end{equation}
Let
$$
V : =\{X=(\vv x_1,\dots,\vv x_m) \  :\   x_{j,i}=0 \  \ \forall\
j={1,\dots,m}\; ; i={2,\dots,n},\ \}  \ ,
$$
where $\vv x_j=(x_{j,1},\dots,x_{j,n})$. Thus, $V$ is an
$m$--dimensional subspace of $\R^{nm}$ and we easily verify
conditions (i) and (ii) of the \MTP{} for linear forms, which is
now applied with $k = m n$,  $l= m(n-1)$ and $B = \I^{nm}$.
Let
$$
\tilde\Psi(\vv q):=g\left(\frac{\Psi_j(\vv q)}{\sqrt{{\vv q}.{\vv
q}}}\right)^{1/m}\sqrt{\vv q.\vv q}\,.
$$
Then,
\begin{equation}\label{e:022}
    \cA_{n,m}^{\vv b}(\tilde\Psi) \ = \ \La(g(\Upsilon)^{1/m})  \, \cap \,
\I^{nm} \ \ .
\end{equation}
Since $r^{-m}g(r)$ is decreasing we have that $$ \tilde\Psi(\vv
q)^m\asymp g\left(\frac{\Psi_j(\vv q)}{|\vv q|}\right)|\vv
q|^m $$ so that $\sum_{\vv q\in\Z^n\nz}\tilde\Psi(\vv q)^m=\infty$.
Therefore, by Sprindzuk's theorem, we have that the set
(\ref{e:022}) has full measure in $\I^{nm}$ and (\ref{e:018}) is
fulfilled. Thus we have (\ref{e:019}). The inclusion
(\ref{e:021}) completes the argument.

As a consequence of Theorem~\ref{bv100} we have the following
statement for approximating functions with sup norm argument. The
case $n=1$, not covered by Theorem~\ref{bv100},
corresponds to Theorem~\ref{galyest}.

\begin{theorem}
\label{mineyourstheirs} Let $n+m>2$. Let $f$ and $g : r \to g(r):=
r^{-m(n-1)} \, f(r)$ be dimension functions such that
$r^{-mn}f(r)$ is monotonic. Let $\psi \, : \, \N \to \Rp $ be an
approximating function. If $n=2$, assume  that  $\psi$ is
monotonic. Then
$$
\cH^f(\cA_{n,m}(\psi))=\left\{\begin{array}{ll} 0 & \ds\text{if }
\;\;\; \sum_{r=1}^\infty \;  g\left(\frac{\psi(r)}{r}\right) \
r^{n+m-1} <\infty \, , \\[2ex] &
\\ \cH^f(\I^{nm}) & \ds\text{if } \;\;\; \sum_{r=1}^\infty \;
g\left(\frac{\psi(r)}{r}\right) \  r^{n+m-1} =\infty\,.
\end{array}\right.
$$
\end{theorem}

We note that the validity of Conjecture \ref{B} together with \MTP{}
for linear forms would remove the monotonicity condition on $\psi$
in the above theorem. With $\psi$ monotonic, the theorem corresponds
to the linear forms  analogue of the Khintchine-\Jarnik{} theorem as
first established in \cite{BDV06}. Obviously, this can be deduced
directly from Groshev's theorem.

Finally, it is easily verified that the \MTP{} for linear forms yields
the following generalisations of the Duffin-Schaeffer and Catlin conjectures
stated in \S\ref{Sprindzuklf}. In short, the Lebesgue conjectures imply the
corresponding  Hausdorff conjectures.

\begin{conjecture}[General Duffin-Schaeffer]
Let $f$ and $g : r \to g(r):=
r^{-m(n-1)} \, f(r)$ be dimension functions such that
$r^{-mn}f(r)$ is monotonic. Let $\Psi:\Z^n\to\Rp$ be an
approximating function. Then
$$
\cH^f(\cA''_{n,m}(\Psi))= \cH^f(\I^{nm})\qquad \ds\text{if }
\;\;\; \sum_{\vv q\in\Z^n\nz} \; g\left(\frac{\Psi(\vv q)}{|\vv
q|}\right)\times \left(\frac{\varphi(\gcd(\vv q))}{\gcd(\vv
q)}|\vv q|\right)^m =\infty\,.
$$
\end{conjecture}

\begin{conjecture}[General Catlin]
Let $f$ and $g : r \to g(r):=
r^{-m(n-1)} \, f(r)$ be dimension functions such that
$r^{-mn}f(r)$ is monotonic. Let $\Psi:\Z^n\to\Rp$ be an
approximating function. Then
$$
\cH^f(\cA'_{n,m}(\Psi))= \cH^f(\I^{nm})\qquad \ds\text{if } \;\;\;
\sum_{\vv q\in\Z^n\nz} \; \max_{t\in\N}g\left(\frac{\Psi(t\vv
q)}{t|\vv q|}\right)\times \cN_{n,m}(\vv q) =\infty\  ,
$$
where $\cN_{n,m}$ is defined in (\ref{e:015}).
\end{conjecture}

\section{Twisted `inhomogeneous' approximation \label{inhomapprox} }

Throughout this section $\psi: \N \to \R^+ $ will be a {\bf
monotonic} approximating function and we write $\cA_{n,m}^{\vv
b}(\psi)$  for the general `inhomogeneous' set $\cA_{n,m}^{\vv
b}(\Psi)$; i.e. $\cA_{n,m}^{\vv b}(\psi)$  is given by
(\ref{e:008}) with $\Psi(\vv q)=\psi(|\vv q|)$. The following
clear cut statement, which is a direct consequence of the
discussion above, provides a complete metric theory for
$\cA_{n,m}^{\vv b}(\psi)$.

\begin{theorem}
\label{inhomreg}
 Let $f$
and $g : r \to g(r):= r^{-m(n-1)} \, f(r)$ be dimension functions
such that $r^{-mn}f(r)$ is monotonic. Let $\psi$ be a monotonic approximating
function. Then
$$
    \cH^f(\cA_{n,m}^{\vv b}(\psi)) = \begin{cases} 0
      &\ds\text{if } \quad \sum_{r = 1}^{ \infty } \
      g\left(\frac{\psi(r)}{r}\right)\times r^{n+m-1}<\infty \ , \\[4ex]
      \cH^f(\I^{nm})
      &\ds\text{if }  \quad \sum_{r = 1}^{ \infty } \
      g\left(\frac{\psi(r)}{r}\right)\times r^{n+m-1} =\infty \ .
                  \end{cases}
$$
\end{theorem}

In view of the discussion in \S \ref{MTPSLF}, this general
Hausdorff measure statement is easily seen to be a consequence of
the corresponding Lebesgue statement (i.e. $f(r) = r^{mn}$ in
Theorem \ref{inhomreg}) and the Mass Transference Principle for
systems of linear forms. It is also worth mentioning, especially
in the context of what is about to follow, that the behavior of
$\cH^f(\cA_{n,m}^{\vv b}(\psi))$ is completely independent of the
fixed inhomogeneous factor  $ {\vv b} \in \I^m$.

We now consider a somewhat `twisted' version of the set
$\cA_{n,m}^{\vv b}(\psi)$ in which the inhomogeneous factor  $
{\vv b}$ becomes the object of approximation. More precisely,
given $ X \in \I^{mn} $ let
$$
\cV^X_{n,m}(\psi) := \{{\vv b} \in \I^m:\|\q X +\vv b \| <
\psi(|\q|) {\text { for
    infinitely many }} \q \in \Z^n \} \ .
$$

\noindent For the ease of motivation and indeed clarity of results
we begin by describing the one--dimensional situation.

\subsection{The one--dimensional theory
\label{oneinhomapprox} }

For any irrational $x$ and any real number $b$, a theorem of
Khintchine states that there are infinitely many integers $q$ such
that
\begin{equation}\label{e:023}
\|q x - b \| < \frac{1+\epsilon}{\sqrt5\, q }  \ .
\end{equation}
In this statement $\epsilon >0 $ is arbitrary and apart from this
term it is equivalent to Hurwitz's homogeneous ($b =0$) theorem. A
weaker form, with $3/q$ appearing on the right hand side of
(\ref{e:023}), had been established earlier by Tchebychef. This
enabled him to conclude that for any irrational $x$ the sequence
$\{qx\}_{q \in \N} $ modulo one is dense in the unit interval. Later
this sequence has been shown to be uniformly distributed. In view of
this density result, it is natural to consider the problem of
approximating points in the unit interval with a pre-described rate
of approximation by the sequence $qx\bmod1$. That is to say, to
investigate the set $\cV^x(\psi) := \cV^x_{1,1}(\psi)$.

Before describing a complete metric theory for $\cV^x(\psi)$ we
state two results which are simple consequences of
(\ref{e:023}) and the Mass Transference Principle (Theorem
\ref{mtp0}). Given $\tau\ge 0$, let $ \psi : r \to r^{-\tau}$ and write  $
\cV^x(\tau) $ for $\cV^x(\psi)$.

\begin{theorem}
\label{inhom1}
 Let $x$ be irrational. For   $\tau \ge 1 $,
$$
\cH^{\frac{1}{\tau}} ( \cV^x(\tau)) \  =  \ \cH^{\frac{1}{\tau}} (
\I ) \ \ .
$$
\end{theorem}
It follows directly from the definition of Hausdorff dimension
that $ \dim \cV^x(\tau) \ge 1/\tau $. The complementary upper
bound result is easily establish and as a corollary we obtain the
following statement.

\begin{corollary}[Bugeaud \cite{Bug}  and Schmeling $\&$
Troubetzkoy \cite{Schmeling}] \label{inhom1cor}
 Let $x$ be irrational. For   $\tau > 1 $,
$$
\dim \cV^x(\tau) =  \frac{1}{\tau}  \  \  .
$$
\end{corollary}

\noindent Thus the corollary is a simple consequences of
(\ref{e:023}) and the Mass Transference Principle. Moreover, we
are able to deduce that the Hausdorff measure at the critical
exponent is infinity; i.e. for   $\tau > 1 $, $$
\cH^{\frac{1}{\tau}} ( \cV^x(\tau)) \ = \infty  \ \ .
$$

\noindent Next, for  $\epsilon >0 $  let $\psi_\epsilon :  r \to
\epsilon/r $ and $f_\epsilon :  r \to r / \epsilon $. Note that
$\cH^{f_\epsilon}$ is simply one dimensional Lebesgue measure
scaled by a multiplying factor $1/ \epsilon$. Now on combining
(\ref{e:023}) and the Mass Transference Principle in the
obvious manner, we obtain that for any irrational $x$
$$
\cH^{f_\epsilon} ( \cV^x(\psi_\epsilon) )  \ = \ \cH^{f_\epsilon}
( \I )  \ = \ 1/ \epsilon \ \ .
$$
Alternatively,
$$
\cH^{f_\epsilon} ( \I \setminus \cV^x(\psi_\epsilon) )  \ = \ 0
$$
and via a simple covering argument we deduce that $ \cH^{1} ( \I
\setminus \cV^x(\psi_\epsilon) )  \ = \ 0 \, $. Thus, for any
$\epsilon >0 $ and any irrational $x$ we have that
\begin{equation}\label{e:024}
| \cV^x(\psi_\epsilon) | \ = \ 1 \ \ ,
\end{equation}
and we have given a short and `direct' proof of the following
statement.
\begin{theorem}[Kim \cite{Kim}]
\label{inhom2}
 Let $x$ be irrational. For almost every $b \in \I$
$$
\liminf_{q \to \infty } q || q x - b || \ = \ 0 \ \ .
$$
\end{theorem}

\noindent Note the above theorems and  corollary are statements
for any irrational $x$ -- the relevance  of this will soon become
apparent.

We now turn our attention towards developing a complete metric
theory for $\cV^x(\psi)$. To begin with, let us concentrate on the
Lebesgue theory. On exploiting the limsup nature of the set
$\cV^x(\psi)$, it is easily verified that for any irrational $x$ and
any approximating function $\psi$
$$
| \cV^x(\psi) | \ = \ 0 \quad {\rm if} \quad \ \sum_{r = 1}^{ \infty
} \ \psi(r)   < \infty \  .
$$
It is worth stressing that the choice of the irrational $x$ and
the `convergent' $\psi$ are completely irrelevant.  Naturally, one
may suspect or even expect that $ | \cV^x(\psi) | = 1 $ if the
above sum diverges -- irrespective of the irrational $x$ and the
`divergent' $\psi$. This is certainly the situation in the
`standard' inhomogeneous setup; i.e. for the set $\cA^{b}(\psi)$
the choice of the inhomogeneous factor $b$ and the `divergent'
$\psi$ are completely irrelevant. However, the `twisted' setup
throws up a few surprises which to some extent lead to a `richer'
theory -- particularly in higher dimensions. Regarding the latter,
it is slightly out of place to give a discussion here and we refer
the reader to \cite{Bugch}. The following statement and indeed its
higher dimension analogue (see \S \ref{mninhomapprox}) is due to
Kurzweil \cite{Kur}.

\begin{theorem}[Kurzweil]
\label{1Kurzweil} Let $\psi$ be given. Then for almost all
irrational $x$
$$ | \cV^x(\psi) | \ = \ 1 \quad
\ds\text{if } \quad \ \sum_{r = 1}^{ \infty } \ \psi(r)   = \infty \
 .
$$
\end{theorem}
As already mentioned the complementary convergent part is valid for
all irrational $x$. The above theorem indicates that the set of
irrational $x$ for which we obtain the full measure statement is
dependent on the choice of the `divergent' $\psi$. To clarify this
and to take the discussion further, let $$\cD := \{ \psi :
\textstyle{ \sum_{r = 1}^{ \infty } \ \psi(r)   = \infty } \}   \ .
$$ Thus the set $\cD$ is the set of `divergent' approximating
functions $\psi$. Also, for $\psi \in \cD $ let
$$
\cV(\psi) \ := \ \{ x \in \I : | \cV^x(\psi) | \, = \, 1 \} \ .
$$
With this notation in mind, Theorem \ref{1Kurzweil} simply states
that $ | \cV(\psi) | \, = \, 1 $. Furthermore, Kurzweil in
\cite{Kur} solves a problem of H. Steinhaus by establishing the
following elegant result which  characterizes the set  $\SBA$ of
badly approximable numbers in terms of `twisted' inhomogeneous
approximation.

\begin{theorem}[Kurzweil]
\label{badKurzweil}
$$
\bigcap_{\psi \in \cD} \cV(\psi) \ = \ \SBA \  .
$$
\end{theorem}

\noindent Thus, for any given  irrational $x$ which is not badly
approximable there is a `divergent' approximating function $\psi$
for which $ | \cV^x(\psi) | \, \neq \, 1 $. In other words,
Theorem \ref{1Kurzweil} is in general false for all irrational
$x$. 

The truth of the statement of Theorem~\ref{badKurzweil} can to
some extent be explained by the fact that the distribution of
$qx\bmod1$ is best possible if $x\in \SBA$. More precisely, it is
well known that the discrepancy $D(N)$ of $qx\bmod1$ satisfies
$D(N)\ll\log N$ if $x\in \SBA$ and that $D(N)\gg\log N$ infinitely
often for any real number $x$. A theorem of Schmidt
\cite{schirreg},  building on the pioneering work of Roth
\cite{rothirreg}, states that the latter is indeed the case for
any sequence $x_n\bmod1$.

Before moving onto the general Hausdorff theory, we point out
that Theorem \ref{1Kurzweil} implies  (\ref{e:024}) for almost
almost all $x$. This leads to the weaker `almost all irrational'
(rather than all irrational) version of Theorem \ref{inhom2}. The
point is that for the particular function $\psi_{\epsilon} $
appearing (\ref{e:024}), it is possible to replace `almost all
irrational' by `all irrational' in the statement of Theorem
\ref{1Kurzweil}. This then implies Theorem \ref{inhom2} as stated.

On applying the Mass Transference Principle in the obvious manner,
we are able to deduce from Theorem \ref{1Kurzweil} the following
general Hausdorff measure statement -- the convergent part is
straightforward and is valid for all irrationals.

\begin{theorem}
\label{1Kurzweilgen}  Let $f$  be dimension function such that
$r^{-1}f(r)$ is monotonic and let $\psi$ be a monotonic approximating
function. Then for almost all irrational $x$
$$
\cH^f(\cV^x(\psi)) = \begin{cases} 0
      &\ds\text{if } \quad \sum_{r = 1}^{ \infty } \
      f(\psi(r) ) <\infty \ , \\[4ex]
      \cH^f(\I)
      &\ds\text{if }  \quad \sum_{r = 1}^{ \infty } \
      f(\psi(r)) =\infty \ .
                  \end{cases}
$$
\end{theorem}

It is worth pointing out that in the case $ \psi : r \to r^{-\tau}$
and $f : r \to r^s$, it is possible to strengthen the theorem to all
irrational $x$ -- see Theorem \ref{inhom1}. The key is that in this
situation one can apply the Mass Transference Principle to
(\ref{e:023}), which is valid for any irrational $x$.
This is the reason why we are able to prove a dimension result
(Corollary \ref{inhom1cor}) for any irrational rather than just
almost all irrational. The latter is all that we can obtain from
Theorem \ref{1Kurzweilgen}.

We end our discussion of the one-dimensional `twisted' theory by
attempting to generalize Theorem \ref{badKurzweil}.  Let $f$  be a
dimension function such that $r^{-1}f(r)$ is monotonic and let
$$\cD^f := \{ \psi : \textstyle{ \sum_{r = 1}^{ \infty } \
f(\psi(r))   = \infty } \}   \ . $$ The set $\cD^f$ is the  set of
`$f$-divergent' approximating functions $\psi$. Also, for $\psi \in
\cD^f $ let
$$
\cV^f(\psi) \ := \ \{ x \in \I : \cH^f(\cV^x(\psi)) \, = \, \cH^f(\I) \} \ .
$$
In this notation, the divergent part of Theorem \ref{1Kurzweilgen} simply states
that $ | \cV^f(\psi) | \, = \, 1 $. Furthermore, on combining the Mass Transference Principle and
Theorem \ref{badKurzweil} we obtain the following result.

\begin{theorem}\label{more}
Let $f$  be dimension function such that
$r^{-1}f(r)$ is monotonic. Then
$$
    \bigcap_{\psi \in \cD^f} \cV^f(\psi) \ \supseteq \ \SBA \  .
$$
\end{theorem}

The following conjecture is a natural refinement of Theorem~\ref{badKurzweil}.

\begin{conjecture}
Let $f$  be dimension function such that
$r^{-1}f(r)$ is monotonic. Then
$$
\bigcap_{\psi \in \cD^f} \cV^f(\psi) \ = \ \SBA \ .
$$
\end{conjecture}

Note that on combining Theorems~\ref{badKurzweil} and \ref{more}
we obtain the following statement.
\begin{corollary}
\label{conjeasy} Let $\cF$ be the set of dimension functions $f$
such that $r^{-1}f(r)$ is monotonic. Then
$$
\bigcap_{f \in \cF} \bigcap_{\psi \in \cD^f} \cV^f(\psi) = \ \SBA
\ .
$$
\end{corollary}

\subsection{The higher dimensional theory
\label{mninhomapprox} }

Starting with the Lebesgue theory, Kurzweil established the higher dimensional
analogue of Theorem \ref{1Kurzweil}.

\begin{theorem}[Kurzweil]
\label{2Kurzweil} Let $\psi$ be given. Then for almost all
$X \in \I^{mn}$
$$ | \cV^X_{n,m}(\psi) |_{m} \ = \ 1 \quad
\ds\text{if } \quad \ \textstyle{\sum_{r = 1}^{ \infty } \ r^{n-1} \ \psi(r)^m   =
\infty \ } .
$$
\end{theorem}

\noindent On applying the Mass Transference Principle (Theorem \ref{mtp0}),
we obtain the following general Hausdorff statement. The convergence part is again
straightforward and is valid for all $X \in \I^{mn}$.

\begin{theorem}
\label{2Kurzweilgen}  Let $f$  be a dimension function such that
$r^{-m}f(r)$ is monotonic and let $\psi$ be a monotonic
approximating function. Then for almost all $X \in \I^{nm}$
$$
\cH^f(\cV^X_{n,m}(\psi)) = \begin{cases} 0
      &\ds\text{if } \quad \sum_{r = 1}^{ \infty } \
      f(\psi(r) ) \ r^{n-1} <\infty \ , \\[4ex]
      \cH^f(\I^m)
      &\ds\text{if }  \quad \sum_{r = 1}^{ \infty } \
      f(\psi(r)) \ r^{n-1}  =\infty \ .
                  \end{cases}
$$
\end{theorem}

Let $ \psi : r \to r^{-\tau}$ and write
$\cV^X_{n,m}(\tau) $ for $\cV^X_{n,m}(\psi)$. As a consequence of the above theorem we have the
following corollary.

\begin{corollary} \label{inhom2cor}
Let  $\tau > n/m $. Then for almost all $X \in \I^{nm}$
$$
\dim
\cV^X_{n,m}(\tau) =  \frac{n}{\tau}  \   \quad       {\it \ and \ moreover
\ } \quad \cH^{\frac{n}{\tau}} ( \cV^X_{n,m}(\tau) ) \ = \infty  \ \
.
$$
\end{corollary}

The above dimension statement is not new -- see  \cite{Bugch}.  It
is worth stressing that in higher dimensions it is not possible to
obtain  a dimension statement for all `irrational' $X$ as in the
one dimensional theory -- see \cite{Bugch}. The point is that the
higher dimensional analogue of (\ref{e:023}) is only valid for
almost all $X \in \I^{mn}$ rather than all `irrational'  $X \in
\I^{mn}$. Finally, we mention that Kurzweil also obtained the
higher dimensional analogue of Theorem \ref{badKurzweil} and
therefore the analogues of Theorem \ref{more} and Corollary
\ref{conjeasy} for arbitrary $n$ and $m$ are  also possible.

\subsection{Back to algebraic irrationals and Roth again}

We end up this paper with another discussion on the interactions
of Roth's theorem and the metrical theory of Diophantine
approximation. As mentioned in
\S\ref{sec:Dirichl-Roth-metrica-theory}, quadratic real algebraic
numbers are badly approximable and Roth's theorem states that
algebraic numbers of degree $n\ge3$ denoted by $\A_n$ are
relatively badly approximable.
It is also believed that $\A_n$
does not contain badly approximable numbers for $n\ge3$. If the
latter is the case then for any algebraic number $\alpha$ of
degree $\ge3$ we have $\alpha\not\in\cV(\psi)$ -- see
Theorem~\ref{badKurzweil}. In other words, one must be able to
construct a monotonic approximating function $\psi$ with
$\sum_{r=1}^\infty\psi(r)=\infty$ such that
$$
|\cV^\alpha(\psi)|<1\,.
$$
Restating this by making use of the definition of the set
$\cV^\alpha(\psi)$ we are naturally led to the following
conjecture.

\begin{conjecture}\label{I}
For any $n\ge3$ and any $\alpha\in\A_n$ there is a monotonic
approximating function $\psi$ and a subset $B$ in $[0,1]$ of
positive Lebesgue measure such that
$$\sum_{r=1}^\infty\psi(r)=\infty$$ but for any $b\in B$ the
inequality
$$
\|q\alpha+b\|<\psi(q)
$$
has only finitely many solutions $q\in\N$.
\end{conjecture}

It is quite possible that one can prove the following
inhomogeneous strengthening of Lang's conjecture that would imply
Conjecture~\ref{I}.

\begin{conjecture}\label{J2}
For any $n\ge3$ and any $\alpha\in\A_n$. Then there is a subset $B$
in $[0,1]$ of positive measure such that for any $b\in B$ the
inequality
$$
\|q\alpha+b\|<\dfrac{1}{q\log q}
$$
has only finitely many solutions $q\in\N$.
\end{conjecture}

Possibly, the only condition that needs to be imposed on $b$ to
fulfil the above conjecture is that $b$ and $\alpha$ are linearly
(or algebraically) independent over $\Q$.

\vspace{4ex}

\noindent\textbf{Acknowledgements.}    SV would like to start by
thanking Klaus Roth for teaching him an important lesson which can
be easily summed up: always think before you write! This was done
in such a nice way that I have never forgotten it. Also his
brilliant talks centred around the theme that there are no
integers between zero and one have been embedded in my grey matter
forever. This latter fact leads me naturally onto Ayesha and Iona
-- so girls, the Scarecrow has finally found his brain and it's
all thanks to Klaus Roth. Thus all is well in Oz and we can all go
home -- whoops this will inevitably lead to some other unknown
journey!

{\small

\providecommand{\bysame}{\leavevmode\hbox to3em{\hrulefill}\thinspace}

\vspace{10mm}

\noindent Victor V. Beresnevich: Department of Mathematics,
University of York,

\vspace{-3mm}

 ~ \hspace{30mm}  Heslington, York, YO10 5DD, England.


 ~ \hspace{30mm} e-mail: vb8@york.ac.uk

\vspace{5mm}

\noindent Vasily I. Bernik:  Institute of Mathematics, Academy of
Sciences of Belarus,

\vspace{-3mm}

  ~ \hspace{19mm}
220072, Surganova 11, Minsk, Belarus.



\vspace{5mm}

\noindent Maurice M. Dodson: Department of Mathematics, University
of York,

\vspace{-3mm}

 ~ \hspace{27mm}  Heslington, York, YO10 5DD, England.


 ~ \hspace{27mm} e-mail: mmd1@york.ac.uk

\vspace{5mm}

\noindent Sanju L. Velani: Department of Mathematics, University
of York,

\vspace{-3mm}

 ~ \hspace{19mm}  Heslington, York, YO10 5DD, England.


 ~ \hspace{19mm} e-mail: slv3@york.ac.uk

}


\begin{thebibliography}{10}

\bibitem{BS}
A.~Baker and W.~M. Schmidt:  \emph{Diophantine approximation and
{H}ausdorff
  dimension}, Proc.\ Lond.\ Math.\ Soc. \textbf{21} (1970), 1--11.

\bibitem{BDV06}
V.~Beresnevich, D.~Dickinson and S.~L. Velani: \emph{Measure
{T}heoretic
  {L}aws for limsup {S}ets}, Mem. \ Amer. \ Math. \ Soc. \textbf{179} (2006),
  no.~846, 1--91.

\bibitem{BV06}
V.~Beresnevich and S.~L. Velani: \emph{A {M}ass {T}ransference
{P}rinciple and
  the {D}uffin-{S}chaeffer conjecture for {H}ausdorff measures}, Ann. \ Math.
  \textbf{164} (2006), 971--992.

\bibitem{BV06Slicing}
\bysame \  \emph{Schmidt's theorem, Hausdorff Measures and
Slicing}, IMRN 2006, Article ID 48794, 24 pages.


\bibitem{BVparis}
\bysame \ \emph{Ubiquity and a general logarithm law for
geodesics}, Pre-print (2007), 22 pages.



\bibitem{Bes34}
A.~S. Besicovitch: \emph{Sets of fractional dimensions ({I}{V}):
on rational
  approximation to real numbers}, J.\ Lond.\ Math.\ Soc. \textbf{9} (1934),
  126--131.

\bibitem{Borel12}
E.~Borel: \emph{Sur un probl{\`e}me de probabilit{\'e}s aux
fractions
  continues}, Math.\ Ann. \textbf{72} (1912), 578--584.



\bibitem{Bug}
Y.~Bugeaud: \emph{A note on inhomogeneous Diophantine
approximation}, Glasgow Math.\ J. \textbf{45} (2003), 105--110.


\bibitem{Bugch}
Y.~Bugeaud and N.~Chevallier:  \emph{ On simultaneous
inhomogeneous Diophantine approximation},   Acta Arith.
\textbf{123} (2006), 97-123.

\bibitem{Casselshort}
J.~W.~S. Cassels: \emph{An introduction to {D}iophantine
approximation},
  Cambridge University Press, 1957.

\bibitem{catlin76}
P.~A. Catlin: \emph{Two problems in metric {D}iophantine
approximation {I}},
  J.~Number Th. \textbf{8} (1976), 282--288.

\bibitem{DavSchmidt70a}
H.~Davenport and W.~M. Schmidt: \emph{Dirichlet's theorem on
{D}iophantine
  approximation}, Inst.\ Alt.\ Mat.\ Symp.\ Math. \textbf{4} (1970), 113--132.

\bibitem{HDSV97}
H.~Dickinson and S.~L. Velani: \emph{Hausdorff measure and linear
forms}, J.\
  reine angew.\ Math. \textbf{490} (1997), 1--36.

\bibitem{DuffinSchaeffer}
R.~J. Duffin and A.~C. Schaeffer: \emph{Khintchine's problem in
metric
  {D}iophantine approximation}, Duke Math.\ J. \textbf{8} (1941), 243--255.

\bibitem{FalcGFS}
K.~Falconer: \emph{The geometry of fractal sets}, Cambridge Tracts
in
  Mathematics, No.~85, Cambridge University Press, 1985.

\bibitem{Gallagher65}
P.~X. Gallagher: \emph{Metric simultaneous {D}iophantine
approximation {II}},
  Mathematika \textbf{12} (1965), 123--127.

\bibitem{HarmanMNT}
G.~Harman: \emph{Metric number theory}, LMS Monographs New Series,
vol.~18,
  Clarendon Press, 1998.

\bibitem{Ja29}
V.~Jarn{\'\i}k: \emph{Diophantischen {A}pproximationen und
{H}ausdorff\-sches
  {M}ass}, Mat.\ Sbornik \textbf{36} (1929), 371--382.

\bibitem{Ja31}
\bysame \  \emph{{\"U}ber die simultanen diophantischen
{A}pproximat\-ionen},
  Math.\ Z. \textbf{33} (1931), 505--543.



\bibitem{Kh24}
A.~I. Khintchine: \emph{Einige {S}{\"a}tze {\"u}ber
{K}ettenbruche, mit
  {A}nwendungen auf die {T}heorie der {D}iophantischen {A}pproximationen},
  Math.\ Ann. \textbf{92} (1924), 115--125.

\bibitem{Kh25}
\bysame \  \emph{{\"U}ber die angen{\"a}herte {A}ufl{\"o}sung
linearer
  {G}leichungen in ganzen {Z}ahlen}, Rec.\ math.\ Soc.\ Moscou Bd. \textbf{32}
  (1925), 203--218.

\bibitem{Kim}
D.~Kim: \emph{The shrinking target property of irrational
rotations}, Pre-print (2007), 1--9.


\bibitem{Kur}
J.~Kurzweil: \emph{On the metric theory of inhomogeneous
Diophantine approximations}, Studia mathematica  \textbf{XV}
(1955), 84-112.



\bibitem{Lang}
S.~Lang: \emph{Report on diophantine approximations}, Bull.\ de la
Soc.\ Math.\ de France 93, (1965), 117--192.


\bibitem{MattilaGS}
P.~Mattila: \emph{Geometry of sets and measures in {E}uclidean
space},
  Cambridge University Press, 1995.

\bibitem{PV90}
A.~D. Pollington and R.~C. Vaughan: \emph{The $k$-dimensional
{D}uffin and
  {S}chaeffer conjecture}, Mathematika \textbf{37} (1990), 190--200.



\bibitem{rothirreg}
K.~F. Roth: \emph{On irregularities of distribution.},
  Mathematika  \textbf{7} (1954), 73--79



\bibitem{Roth}
\bysame \   \emph{Rational approximation to algebraic numbers},
Mathematika
  \textbf{2} (1955), 1--20, with a corrigendum on p.~168.


\bibitem{Schmeling}
J.~Schmeling and S.~Troubetzkoy: \emph{Inhomogeneous Diophantine
Approximation and Angular Recurrence for Polygonal Billards},
Math.\ Sbornik \textbf{194} (2003), 295--309.

\bibitem{Schmidt64}
W.~M. Schmidt: \emph{Metrical theorems on fractional parts of
sequences},
  Trans.\ Amer.\ Math.\ Soc. \textbf{110} (1964), 493--518.


\bibitem{schirreg}
\bysame \  \emph{Irregularities of Distribution. VII.},
  Acta.\ Arith. \textbf{21} (1972), 45---50.


\bibitem{Schmidtbook}
\bysame \  \emph{Diophantine approximation}, Lecture notes in
sMath. 785, Springer -- Verlag, (1975).


\bibitem{Sprindzuk}
V.~G. Sprind{\v z}uk: \emph{Metric theory of {D}iophantine
approximations},
  John Wiley, 1979, Translated by R.~A.~Silverman.

\bibitem{Waldschmidt}
M.~Waldschmidt: \emph{Open Diophantine Problems}, Moscow
Mathematical Journal \textbf{4}  (2004), 245–305.

\end{thebibliography}
\end{document}